\documentclass{article}%
\usepackage{amsfonts}
\usepackage{amsmath}
\usepackage{amssymb}
\usepackage{graphicx}%
\usepackage{cite}
\usepackage[usenames]{color}
\usepackage[normalem]{ulem}
\providecommand{\U}[1]{\protect\rule{.1in}{.1in}}

\usepackage{orcidlink}

\begin{document}

\title{Hyperbolic extensions of constrained PDEs}
\author{Fernando Abalos\orcidlink{0000-0001-7863-3711} \footnote{Email: j.abalos@uib.es} $\,{}^{1}$, Oscar Reula\orcidlink{0000-0003-2517-7454}$\,{}^2$ and David Hilditch\orcidlink{0000-0001-9960-5293}$\,{}^3$}
\date{\small{
${}^1$Departament de F\'isica, Universitat de les Illes Balears, Palma de Mallorca, E-07122, Spain and Institute of Applied Computing and Community Code (IAC3),
Universitat de les Illes Balears, Palma de Mallorca, E-07122, Spain,\\
${}^2$Facultad de Matem\'atica, Astronom\'ia y F\'isica, Universidad Nacional de C\'ordoba and IFEG-CONICET, Ciudad Universitaria, X5016LAE C\'ordoba, Argentina,\\
${}^3$CENTRA, Departamento de F\'isica, Instituto Superior T\'ecnico IST,
Universidade de Lisboa UL, Avenida Rovisco Pais 1, 1049 Lisboa, Portugal}
}

\maketitle

\begin{abstract}
Systems of PDEs comprised of a combination of constraints and evolution equations are ubiquitous in physics. For both theoretical and practical reasons, such as numerical integration, it is desirable to have a systematic understanding of the well-posedness of the Cauchy problem for these systems. Presently we review the use of hyperbolic reductions, in which the evolution equations are singled out for consideration. We then examine in greater detail the extensions, in which constraints are evolved as auxiliary variables alongside the original variables. Assuming a particular structure of the original system, we give sufficient conditions for strong-hyperbolicity of an extension. This theory is then applied to the examples of electromagnetism and a toy for magnetohydrodynamics.

\end{abstract}

\section{Introduction}

In this work we continue~\cite{geroch1996partial,abalos2017necessary,Abalos:2018uwg,Abalos_2022} the study of first order systems of equations in which there are more equations than unknowns, but with a structure which permits, in principle, to split suitable linear combinations of them into ``evolution'' and ``constraint'' equations. We restrict to the case of consistent systems, in which the number of equations is equal to the number of constraints plus the number of dependent variables, and furthermore to the special case in which the number of dependent variables matches the number of evolution equations. The latter means that we do not consider systems with gauge freedom remaining, which would imply the existence of variables with unspecified equations of motion. In this case one can attempt to solve by carefully restricting the initial data and then directly solving the evolution equations. For an introductory review, see~\cite{david_13}. One must then check that the constraint equations are satisfied in the time development. For this, integrability identities among the whole system of equations must be satisfied. These conditions will be assumed and spelled out in detail below. This ``free evolution approach" requires us to establish well-posedness of the Cauchy problem~\cite{gustafsson1995time,kreiss2014introduction} (for a review about well-posedness applied to general relativity see~\cite{sarbach2012continuum}). We restrict ourselves to the concepts arising from the theory of strongly hyperbolic systems, in which well-posedness is determined by algebraic properties of the principal symbol of the equation system. For first order systems the principal symbol this is just the set of matrices multiplying the derivatives of the variables. The algebraic properties leading to well-posedness have several equivalent characterizations summarized in the Kreiss matrix theorem~\cite{kreiss1962stabilitatsdefinition}. To assert well-posedness for the systems under consideration we need to find a suitable square system, that is, a system where the number of variables equals the number of equations. This can be achieved by taking a subset of the equation system, called a \textit{reduction}, resulting in a pure evolution system. The use of reductions is customary but another possibility, which is often employed in numerical schemes, consists of making  an \textit{extension}, that is extending the system by adding more variables. These extensions are commonly referred to as \textit{divergence cleaning} \cite{dedner2002hyperbolic, MUNZ2000484, munz2000three} from its use in magnetohydrodynamics, or as \textit{$\lambda$}~\cite{brodbeck1999einstein} or \textit{Z-systems}~\cite{Bona:2002ft} from their use in general relativity.

A paradigmatic example is given by the Maxwell equations,
\begin{align}
    \nabla_a F^{ab} = J^b\,, 
    \qquad \varepsilon^{dabc}\nabla_{a} F_{bc} = 0\,, 
    \qquad \nabla_a J^a = 0\,,
\end{align}
where the unknowns are the components of the Faraday tensor~$F_{ab}$, an anti-symmetric tensor (so there are a total of 6 dependent variables), $J^a$, the current vector is a given vector fixed in spacetime, but which has vanishing divergence. This is necessary due to the integrability identity~$\nabla_b(\nabla_a F^{ab})=0$. We work here in four dimensional spacetime~$(M,g_{ab})$ with Levi-Civita derivative~$\nabla_a$ associated with~$g_{ab}$. There are thus a total of $8=4+4$ equations for~$F^{ab}$, so 6 of them should be evolution equations and the remaining 2 should be constraints. Introducing a timelike covector~$n_a$ one finds that contraction with that vector on both equations gives constraint, that is, equations which have derivatives only in directions perpendicular to~$n_a$; while projection on the space perpendicular to $n_a$ gives equations which have derivatives along~$n^a$ for each of the independent components of $F^{ab}$.  Thus, in the terminology introduced above, a reduction is obtained by taking just these projections as the evolution equations. The integrability identity and divergence property of~$J^a$ together imply that constraints are satisfied in the time development if they are at an initial surface. 

On the other hand, an extension is given by adding 2 auxiliary constraint variables~$(Z_1,Z_2)$, one for each Maxwell constraint, and making a choice for their equations of motion. To accomplish this in a covariant fashion we need to define 2 tensor fields $(g_1, g_2)$. The proposed extended system is
\begin{align}
    \nabla_a F^{ab} + g_1^{ba}\nabla_a Z_1 = J^b\,, 
    \qquad 
    \varepsilon^{dabc}\nabla_{a} F_{bc}  + g_2^{ba}\nabla_a Z_2 = 0\,, 
    \qquad 
    \nabla_a J^a = 0\,, \label{eqn:Max}
\end{align}
It turns out that if the symmetric parts of~$(g_1, g_2)$ are Lorentzian metrics whose cones have non-zero intersections among each other and with the cone of~$g$, then the extended system is well-posed. (We use the mathematical notion of cone; when needed we use the term light cone to refer to their boundaries). The equations that were constraints are now evolution equations for~$(Z_1,Z_2)$, and the others acquire spatial derivatives of these fields. As mentioned above, such extensions have been employed with enormous success in numerical relativity~\cite{Baumgarte:1998te,Shibata:1995we,NakOohKoj87,Pretorius:2004jg,Bernuzzi:2009ex,alic2012conformal} and computational astrophysics, with works introducing this approach for magnetohydrodynamics~\cite{munz2000three,MUNZ2000484,dedner2002hyperbolic} particularly influential. Here we investigate the space of possible extensions that lead to well posed Cauchy problems and how to construct them in a natural, covariant fashion.

The paper is organized as follows. In Section~\ref{Section:Preliminaries_Notation}, we define the type of systems to be considered, including necessary conditions they must satisfy in order to have a well-posed Cauchy problem. In Section~\ref{Section:Kronecker}, we introduce the Kronecker decomposition of matrix pencils and explain its implications to the study of strongly hyperbolic systems. In Section~\ref{Section:Extensions}, we formalize the framework for extensions. Given the considerable freedom in choosing them, we use the Kronecker decomposition as a guide for making these choices. In Section~\ref{Section:Examples}, we demonstrate how this framework applies to two concrete examples: Maxwell's Electrodynamics and a toy model of MHD. Finally, in Section~\ref{Section:Conclusions}, we conclude with discussions and provide comments on how this line of research are being further developed.

\section{Preliminaries and notation}
\label{Section:Preliminaries_Notation}

To fix notation we specify the systems we consider, following the notation of~\cite{geroch1996partial,Abalos:2018uwg,Abalos_2022}. We consider a manifold~$M$ of dimension~$n$, and the following system of the quasi-linear first order partial differential equations on the fields $\phi$,
\begin{align}
E^{A}:=\mathfrak{N}^{Aa}{}_{\alpha}\left(  x,\phi\right)  \nabla_{a}%
\phi^{\alpha}-J^{A}\left(  x,\phi\right)  =0, \label{eq_sys_1}%
\end{align}
where the indices~($A$, $a$, $\alpha$) are abstract, grouping several tensorial indices into one, and merely indicate where the contractions are. We use lower-case Latin indices to denote single vector indices, lower-case Greek indices to indicate variable fields, and upper-case Latin to label the equations space. The $|\cdot |$ function on indices indicates their total dimension.

We impose the following on~$\mathfrak{N}^{Aa}{}_{\alpha}\left(x,\phi\right)$:

\textbf{Condition 1: The Generalized Kreiss Condition.}

We assume that the matrix~$\mathfrak{N}^{Aa}{}_{\alpha}\left(x,\phi\right)$ is smooth in all arguments, and that there exists a hypersurface orthogonal covector~$n_a$ such that for all values of~$k_a$, not proportional to~$n_a$, the matrix pencil
\begin{align*}
\mathfrak{N}^{Aa}{}_{\alpha}l_a(\lambda)=\mathfrak{N}^{Aa}{}_{\alpha}(\lambda n_a + k_a),
\end{align*}
has kernel only for a finite set of real values~$\{\lambda_i(k)\}$ of $\lambda$ (the term matrix pencil refers here to the uni-parametric combination $\lambda \mathfrak{N} + \mathfrak{B}$, where~$\mathfrak{N}$ and~$\mathfrak{B}$ are matrices that do not depend on the parameter~$\lambda$).  

In addition, the corresponding singular values of~$\mathfrak{N}^{Aa}{}_{\alpha}l_a(\lambda)$ approach zero in a linear way, that is, $\sigma(\lambda) \geq c_i|(\lambda - \lambda_i)|,$ with~$c_i>0$ in a neighborhood of~$\lambda_i$. We recall that the singular values are the square roots of the eigenvalues of~$(\mathfrak{N}^{Aa}{}_{\alpha}l_a)^{T} (\mathfrak{N}^{Ab}{}_{\beta}l_b)$, since this is an~$|\alpha| \times |\alpha|$ matrix, there are~$|\alpha|$ singular values (see~\cite{abalos2017necessary} for more details and for a more general definition).

These conditions imply two things: \textit{ i}) the rank of~$\mathfrak{N}_{~\alpha}^{Aa}\left(  x,\phi\right)n_a$ is maximal, therefore defining any vector~$t^a$ transversal to the surface flat defined by~$n_a$ (i.e. $t^a n_a\neq 0$) we can obtain all field derivatives along~$t^a$ from their values and of their derivatives at that surface. This means that we have enough evolution equations for each field~$\phi^{\alpha}$. Observe that once we have a choice of~$n_a$ satisfying condition~1, then there is an open set of covectors satisfying the same condition, thus we can always form hypersurfaces in a neighorhood of any point, leading to a local initial value problem; \textit{ii}) In the case that the number of equations equals the number of variables, these conditions imply there is a well-posed Cauchy problem, in the usual sense for strongly hyperbolic systems, off of the mentioned surface. This is the classic Kreiss condition.

In case there are more equations than variables, we need to make sure that there are no more linearly independent equations having derivatives along the transversal vector~$t^a$, otherwise we would have an inconsistency because two equations could give different values for the same transversal derivative. To accomplish that we impose:

\textbf{Condition 2: The Geroch Constraint Condition.} 

If the number of equations is larger than the number of variables~$|A| > |\alpha|$, then we assume there exist a set of matrices~$C^{\Gamma a}{}_A$, which are labelled by upper case Greek indices, with 
\begin{align*}
C^{\Gamma (a}{}_{A}\mathfrak{N}^{\left\vert A\right\vert b)}{}_{\alpha}=0,
\end{align*}
and that $\text{rank}(C^{\Gamma a}{}_{A}n_a )= |A| - |\alpha|  = |\Gamma|$. 
This condition ensures that the rest of the equations do not have derivatives off of the surface defined by~$n_a$, so that the system is consistent. Indeed the following linear combination of equations, called constraints,
\begin{align*}
 \psi^{\Gamma} := C^{\Gamma a}{}_{A}n_a(\mathfrak{N}^{A  b}{}_{\alpha}\nabla_b \phi^{\alpha} - J^A),
\end{align*}
have only derivatives on the flat defined by $n_a$.

There is a further consistency condition, one that would guarantee that if the initial data is such that constraint quantities vanish at the initial surface, then they would vanish also along evolution \cite{Abalos_2022}. We require:

\textbf{Condition 3: Integrability.} 

\begin{align}
\nabla_{d}\left(  C^{\Gamma d}{}_{A}E^{A}\right)=L_{1A}^{\Gamma}\left(
x,\phi,\nabla\phi\right)  E^{A}\left(  x,\phi,\nabla\phi\right)  ,
\label{eq_int_G_1}%
\end{align}
In other words, there is a particular on-shell identity among derivatives of our equation system. In most cases of physical interest this identity is a consequence of gauge or diffeomorphism invariance. 

\section{Kronecker decomposition}
\label{Section:Kronecker}

When studying well-posedness of the Cauchy problem the relevant aspect is the behaviour of the system in the limit of high frequencies. We can thus restrict our attention to a neighborhood of each point and work in the frequency domain, employing the Fourier-Laplace transform in space and time respectively. Explicitly we consider a time function~$t$ and a foliation given by its level surfaces. We define~$n_a=(dt)_a$  and take a vector~$t^a$ transversal to the foliation and adjust it such that~$t^an_a = 1$, we choose covectors~$k_a$ such that~$t^ak_a = 0$ and define~$l_a = \lambda n_a + k_a$, we perform Fourier in~$k_a$, and Laplace in~$\lambda$. Thus we replace space derivatives by~$ik_a$ and time derivatives by~$i\lambda$. Furthermore in what follows, once any particular~$k_a$ is chosen we take a coordinate base so that~$n_a = (dx^0)_a$, and~$k_a = (dx^1)_a$, and so~$l_a = (\lambda n_a + k_a) = (\lambda dx^0 + dx^1)_a$. Finally in the  high frequency limit  we obtain~$\mathfrak{N}^{Aa}{}_{\alpha}l_a\tilde{\phi}^\alpha=0$.

The Kronecker decomposition of a matrix pencil is a canonical transformation that generalizes the Jordan decomposition of a square matrix pencil. Considering the (square or non-square) pencil~$\mathfrak{N} \lambda + \mathfrak{B}$, the Kronecker decomposition is achieved by multiplying this pencil by specific matrices~$W$ and~$Q$, which are independent of~$\lambda$ (as in the square Jordan decomposition case). This transformation results in a new pencil~$(W \mathfrak{N} Q) \lambda + (W \mathfrak{B} Q)$ that has a block structure with particular canonical blocks (see~\cite{gantmacher1992theory,gantmakher1998theory} for a detailed description and Eq.~\eqref{KD_K} for an example). 

It turns out that the Kronecker decomposition can be used naturally in the analysis of systems with constraints or gauge freedom. With the first two conditions assumed above the Kronecker decomposition of the pencil~$\mathfrak{N}^{Aa}{}_{\alpha}l_a(\lambda)$ is given by
\begin{align}
\mathfrak{N}^{Aa}{}_{\alpha}l_{a}=\left[
\begin{array}
[c]{ccccccc}%
\lambda-\lambda_{1} & 0 & 0 & 0 & ... & &0 \\
0 & ... & 0 &  &  & &0\\
0 & 0 & \lambda-\lambda_{d} &  0 & ... & &0 \\
0 & ... & 0  & \lambda & 0 & ... & 0 \\
0 & ... & 0 & 1 & 0 & ... & 0\\
0 & ... & 0 & 0 & \lambda &... & 0 \\
0 & ... & 0 & 0 & 1 &... & 0 \\
0 & ... & & ... & &... & 0 \\
0 & ... & & ... & &... & 0 \\
0 & ... & &... && 0 & \lambda\\
0 & ... & &...&& 0 &1\\
0 & 0 & 0 & ... & 0 & 0 & 0\\
... & ... & ... & ... & ... & ...\\
0 & 0 & 0 & ... & 0 & 0 & 0
\end{array}
\right] \label{KD_K}
\end{align}

Ultimately this represents a change of basis both of the variable and equation spaces which depends on~$k_a$, but not on~$\lambda$. The first block is a diagonal~$d \times d$ block, they represent the true degrees of freedom of the whole system, there are as many elements as the zeros of the singular value decomposition counting their multiplicity. The~$1 \times 2$ blocks, called~$L_1^{T}$ in the literature, are due to the constraints, there are a total of~$r =  |\alpha| - d$ blocks. Since each block occupy 2 rows, we see that the number of zero rows is:~$s = |A| - d - 2r$. In many systems the zero rows are present, they represent differential constraints among the constraints themselves. The numbers defined above also satisfy: 
\begin{align}
d\left(  v\right)  &:=\dim\left(  \text{right}\_\ker\left(  C^{\Gamma
a}{}_{A}  n_a\mathfrak{N}^{Ai}{}_{\alpha}k_{i}\right)  \right)  ,\label{Eq_d_1}\\
r\left(  v\right)  &:=\text{rank}\left(  C^{\Gamma a}{}_{A}n_a\mathfrak{N}
^{Ai}{}_{\alpha}k_{i}\right)  ,\label{Eq_r_1}\\
s\left(  v\right)  &=:\dim\left(  \text{left}\_\ker\left(  C^{\Gamma a}{}_{A}n_a
\mathfrak{N}^{Ai}{}_{\alpha}k_{i}\right)  \right)  .\label{Eq_s_1}
\end{align}
With this decomposition at hand, it is easy to see how to choose among them linear combinations which give evolution equations for all~$\phi^{\alpha}$. Observe that the equations (rows) with a~$\lambda$ are certain to contain derivatives transversal to the~$n_a$ flats. So we must include them, but we can add to them any combination of the other rows. It turns out that just adding to each one of these rows the immediate one below, multiplied by any number $\{\pi_i\}, \; i=1,\ldots,r$, and throwing all the remaining rows we get the evolution equations
\begin{align*}
h^{\beta}{}_{A}\mathfrak{N}^{Aa}{}_{\alpha}l_{a}:=\left[
\begin{array}
[c]{cccccc}%
\lambda-\lambda_{1} & 0 & 0 &  &  & \\
0 & ... & 0 &  &  & \\
0 & 0 & \lambda-\lambda_{d} &  &  & \\
&  &  & \lambda - \pi_1 & 0 & 0 \\
&  &  &  & ... & \\
&  &  & 0 &  0 & \lambda - \pi_r\\
\end{array}
\right],
\end{align*}
where we have constructed a map from the equation space to the  variable space which we call a reduction and denote by~$h^{\beta}{}_{A}$. Thus, $h^{\beta}{}_{A}\mathfrak{N}^{Aa}{}_{\alpha}l_{a}$ is a map from the variable space into itself consisting on a set of diagonal matrices satisfying the classic Kreiss conditions (see point~\textit{ii}. within Condition~1). Notice that we can choose the extra roots of~$\lambda$ (i.e the $\{\pi_i\}$) as we please. They are the propagation speed of extra constraint modes. This simple observation is the principle behind the results in~\cite{reula2004strongly, Abalos_2022}.

Thus, there is a reduction (a linear combination of the equations) such that the Cauchy problem of the system is well-posed and furthermore \textbf{Condition 3} asserts that if the initial data satisfies all equations (including the vanishing of the constraints) then all the equations are satisfied for all times as long as the solution exists. See~\cite{Abalos_2022} for details.

\section{Extensions}
\label{Section:Extensions}

A generic extension would imply the addition of an extra matrix, $\tilde{\mathfrak{N}}^{\Delta Aa}\left(  x,\phi\right)$ (and extra variables $Z_{\Gamma}$), to obtain a square system
\begin{align}
\mathfrak{N}^{Aa}{}_{\alpha}\left(  x,\phi\right)  \nabla_{a}\phi^{\alpha
}
+\tilde{\mathfrak{N}}^{\Gamma Aa}\left(  x,\phi\right) \nabla_{a}Z_{\Gamma}-J^{A}\left(  x,\phi\right)
+B^{A}\left(  x,\phi,Z\right)  =0\,. \label{ext_1}
\end{align}
Here~$B^{A}\left(  x,\phi,Z\right)$ is a term we can also freely choose which does not include derivatives of $\phi$ or $Z$ and which goes to zero when Z goes to zero. In general, $B^{A}$ represents damping terms~\cite{brodbeck1999einstein,dedner2002hyperbolic,gundlach2005constraint} which are important in numerical applications. For simplicity in our discussion however, we omit it.

Since we are interested in solving Eq.~\eqref{eq_sys_1} for $\phi$,  our extension proposal only makes sense if we can show that for suitable initial data (for $(\phi,Z) $), the solution of Eq.~\eqref{ext_1} has $Z=0$ in the development, and so $\phi$ a solution of (\ref{eq_sys_1}).  

As we explained before, if we assume conditions 1, 2 and 3 and take any initial data for~$\phi$ satisfying the constraints we know that the initial value problem for Eq.~\eqref{eq_sys_1} is ``well-posed'' and has a unique solution~$\phi_{\text{sol}}$. (Here by well-posed we mean that the map from Cauchy data to solutions is continuous. To establish this one finds a hyperbolic reduction from which we may assert that the reduced system is well-posed for arbitrary initial data. Then one shows that if the initial data satisfies the constraints, then the solutions of the reduced system also satisfy them. Thus they are solutions to the whole system, and we call the whole system well-posed). Therefore if we choose $\tilde{\mathfrak{N}}^{\Gamma Aa}$ such that the extended system, Eq.~\eqref{ext_1}, is well-posed then for any initial data there will be a unique solution. If we choose as initial data $ (\left. \phi_{\text{sol}}\right\vert _{t=0}, \left.  Z\right\vert _{t=0}=0 )$ then $(\phi_{\text{sol}},Z=0)$ will be a solution, and by uniqueness is {\it the} solution. Therefore, we only need to show that system~\eqref{ext_1} satisfies Kreiss's condition.  

\subsection{Strong hyperbolicity of the extensions}

A particularly interesting set of extensions is obtained by noticing the symmetry between the Kronecker decomposition of~$\mathfrak{N}^{Aa}{}_{\alpha}l_a(\lambda)$ and~$(C^{\Delta a}{}_{B}l_a(\lambda))^{T}$. So we start by computing it: 

\begin{align*}
\left(  C^{\Gamma b}{}_{A}l_{b}\right)  ^{T}=\left[
\begin{array}
[c]{cccccc}%
0 & 0 & 0 & 0 & 0 & 0\\
0 & 0 & 0 & 0 & 0 & 0\\
0 & 0 & 0 & 0 & 0 & 0\\
-1 &  &  &  &  & \\
\lambda &  &  &  &  & \\
& ... &  &  &  & \\
&  & -1 &  &  & \\
&  & \lambda &  &  & \\
&  &  & \lambda-\rho_{1} &  & \\
&  &  &  & ... & \\
&  &  &  &  & \lambda-\rho_{s}%
\end{array}
\right]
\end{align*}
Recalling that the matrices~$C^{\Delta a}{}_{B}l_a$ can be thought as a basis, labeled by~$\Delta$, for the the kernel of~$\mathfrak{N}^{Aa}{}_{\alpha}l_a$, it is easy to understand its structure. Here the rows with zeros are~$d$ in number, this is so because the diagonal part of~$\mathfrak{N}^{Aa}{}_{\alpha}l_a$ can not contribute to the kernel. We then have~$r$ blocks~$\left[
\begin{array}
[c]{cc}%
-1 & \lambda
\end{array}
\right]^T$, observing that they have a minus sign on them. This is because they are kernels for the corresponding~$L_1^T$  blocks of~$\mathfrak{N}^{Aa}{}_{\alpha}l_a$. Finally there is a block which is kernel of the zero rows of~$\mathfrak{N}^{Aa}{}_{\alpha}l_a$, this part is completely undetermined, so we have just added a diagonal matrix.

To make more apparent the extension we shall propose, we reorganise the rows of~$\mathfrak{N}%
^{Aa}{}_{\alpha}l_{a}$ and $(C^{\Gamma b}{}_{A} l_{b})^T$ such that

\begin{align}
\mathfrak{N}^{Aa}{}_{\alpha}l_{a}=\left[
\begin{array}
[c]{cc}%
J & 0\\
0 & \lambda I_{r}\\
0 & I_{r}\\
0 & 0
\end{array}
\right],  \text{ \ \ \ }\left(  C^{\Gamma b}{}_{A}l_{b}\right)  ^{T}=\left[
\begin{array}
[c]{cc}%
0 & 0\\
-I_{r} & 0\\
\lambda I_{r} & 0\\
0 & J_{c}%
\end{array}
\right].  \label{eq_K_C_KD_1}%
\end{align}
Here all the matrices are blocks matrices where $J=(\lambda-\lambda_1,..., \lambda-\lambda_d )$ of size $ d  \times d  $, $J_{c}=(\lambda-\rho_1,..., \lambda-\rho_s )$ of size $s \times s$, $I_{r}$ is the identity matrix of size $r \times r$, the zero rows of  $\mathfrak{N}^{Aa}{}_{\alpha}l_{a}$  are of size $s  \times\left\vert \alpha\right\vert $ and the zero rows of $\left(  C^{\Gamma b}{}_{A}l_{b}\right)  ^{T}$ are of $r  \times\left\vert \Gamma\right\vert$.

From this reorganization, it is apparent that a natural choice of~$\tilde{\mathfrak{N}}^{\Gamma Aa}$ is given by
\begin{align}
\tilde{\mathfrak{N}}^{\Gamma Aa}=G^{AB}C^{\Gamma a}{}_{B}, \label{K_mono}%
\end{align}
where $G^{AB}$ now has to be chosen to render the system diagonalizable. This is of course not the most general extension but is a natural and fully covariant proposal for~$\tilde{\mathfrak{N}}^{\Gamma Aa}$. The principal symbol of Eq.~\eqref{ext_1} becomes then, 
\begin{align}
M^{Aa}{}_{D}l_{a}=\left[
\begin{array}
[c]{cc}%
\mathfrak{N}^{Aa}{}_{\alpha} & G^{AB}C^{\Delta a}{}_{B}%
\end{array}
\right]  l_{a}, \label{Eq_M_s_1}%
\end{align}
a $\left\vert A\right\vert \times\left\vert
A\right\vert $ square matrix. 

We now propose a particular expression for $G^{AB}$, namely,
\begin{align}
G^{AB}=\left[
\begin{array}
[c]{cccc}%
I_{d} & 0 & 0 & 0\\
0 & -D^{2} & 0 & 0\\
0 & 0 & I_{r} & 0\\
0 & 0 & 0 & I_{s}%
\end{array}
\right]  ,\label{G_AB_1}%
\end{align}
where~$D=\text{diag}\left(  \pi_{1},...,\pi_{r}\right)$ is of size~$r \times r $ and  $I_{s}$ is the identity matrix of size~$s  \times s$.

Using expressions~\eqref{eq_K_C_KD_1} and~\eqref{G_AB_1},  we conclude
\begin{align}
M^{Aa}{}_{D}l_{a}=\left[
\begin{array}
[c]{cccc}%
J & 0 & 0 & 0\\
0 & \lambda I & D^{2} & 0\\
0 & I & \lambda I & 0\\
0 & 0 & 0 & J_{c}%
\end{array}
\right],
\end{align}
and it is easy to check this matrix is pencil-similar to the following diagonal matrix
\begin{align*}
M^{Aa}{}_{D}l_{a} \sim 
\text{diag} \left(\dots,\lambda-\lambda_{i},\dots, \lambda+\pi_j, \lambda-\pi_j,\dots, \lambda -\rho_k,\dots\right)
\end{align*}
and so it satisfies Kreiss's condition. The extra~$2r$ eigenvalues~$\{\pi_i,-\pi_i\} $, introduced by~$G^{AB}$, come in pairs which means that there are~$r$ new null cones as characteristic. We shall see this in the examples below, where Lorentzian metrics are used to realize these null cones.

\section{Examples}
\label{Section:Examples}

In this section we show several implementation examples of our proposal, showing that it leads to well-posed systems which retain, to a large extent, the covariance of the original theories. In all cases introduction of extra Lorentzian metrics is used to avoid light cone intersections.

\subsection{Maxwell's equations}

We start with the example given in the introduction Eq.~\eqref{eqn:Max}. For them we have, a space of variables~$F^{ab}$ (anti-symmetric tensors), which is~$|\alpha| = 6$ dimensional in a four dimensional space-time of metric~$g_{ab}$. The space of equations is~$|A|=8$, namely two space-time vectors. 
We have (see~\cite{geroch1996partial}), 

\begin{align}
    \mathfrak{N}^{Aa}{}_{\alpha} = \left( \begin{array}{c} 
        \delta^a_{[c} \delta^q_{d]} \\
        \varepsilon^{pa}{}_{bc}
    \end{array} \right) \;\;\;\;\; 
    C^{b\Gamma}{}_{A} = \left( \begin{array}{c} 
        \delta^b{}_q \\
        \delta^b{}_{p}
    \end{array} \right)  \;\;\;\;\; 
    C^{b\Gamma}{}_{A} l_b = \left( \begin{array}{c} 
        l_q \\
        l_p
    \end{array} \right)
\end{align}
Given a time-like $n_a$ we have, 

\begin{align*}
\mathfrak{N}^{Aa}{}_{\alpha}n_a  =
\left( \begin{array}{c}
n_{[c} \delta^q_{d]} \\
\varepsilon^{pa}{}_{bc}n_a
\end{array} \right).
\end{align*}
So it is the map~$F_{ab} \to (E_a,B_a)$, which is of the maximal rank. This system satisfies  \textbf{condition 1}, see \cite{Abalos:2018uwg} for more details.

$C^{b\Gamma}{}_{A} l_b$ is also of maximal rank for any $l$'s \footnote{Here the target space is two copies of $R^4$, and the image is 1-dimensional on each one of them.}. Since the dimension of the image is 2-dimensional we have~$|A| = |\alpha| + |\Gamma|$ and the system is consistent, satisfying \textbf{condition 2}.

We also have, 
\begin{align*}
\nabla_b(C^{b\Gamma}{}_A \mathfrak{N}^{A a}{}_{\alpha}\nabla_a \phi^{\alpha})  = \nabla_b \left( \begin{array}{c}
\delta^a_{[c} \delta^b_{d]} \nabla_a F^{cd} \\
\varepsilon^{ba}{}_{cd}\nabla_a F^{cd}
\end{array} \right)=
\left( \begin{array}{c}
\nabla_bJ^b \\
0
\end{array} \right)= 0
\end{align*}
and so \textbf{Condition 3} is also satisfied.

A suitable reduction is
\begin{align}
    h_{\beta B} = (g_{q[r}t_{s]}, -\frac{3}{2} \varepsilon_{pars}t^a).
\end{align}
This renders the evolution equations symmetric-hyperbolic. As we saw above, a simple extension is obtained introducing two tensors~$(g_1^{pq},g_2^{pq})$ and defining,
\begin{align}
G^{AB} = \left( \begin{array}{cc} 
        g_1^{pq} & 0\\
        0 & g_2^{pq}
    \end{array} \right)
\end{align}
If we take their symmetric parts to be any two Lorentzian metrics, each one of them sharing a common time-like co-vector~$n_a$ with~$g_{ab}$, but not touching their null cones (for brevity we do not consider here such degenerate cases), then the system is strongly hyperbolic and so has a well-posed Cauchy problem. To check this we compute the characteristics of the system and the corresponding eigenvectors and see when we get a complete set, that is a total of 8 eigenvectors.

The characteristic equations are
\begin{align}
    l_b \delta F^{ab} + g_1^{ab}l_b \delta Z_1 &= 0 \\
    \varepsilon^{abcd}l_b \delta F_{cd} + g_2^{ab}l_b \delta Z_2&= 0, 
\end{align}
where we need to solve these equations for~$\lambda$ with $l_a=\lambda n_a+k_a$ and~$n_a,k_a$ given, and for the eigenvectors~$\delta F^{ab}$ and~$\delta Z_{1,2}$. The solutions split in to three cases, first when~$l_a$ is null with respect to~$g^{ab}$ (physical case), then when it is null with respect to~$g_1^{ab}$ or~$g_2^{ab}$ (extended cases) as we explain below.

We already know four of the eigenvectors, namely the physical ones arising from the original system. To recover these, we set~$\delta Z_1=\delta Z_2=0$ and search for the value of~$\delta F_{ab}$. The second equation then implies that~$\delta F_{cd} = 2 l_{[c}A_{d]}$ for some vector~$A_d,$ while the first implies that~$(l_al^a) A^b - (l_aA^a) l^b = 0$ where indices are raised with the space-time metric. Since~$A_a$ can not be proportional to~$l_a$ (otherwise~$\delta F_{cd}$ would vanish) both terms must vanish and so we conclude 
\begin{align}
  g^{ab} l_a l_b = 0\,, \label{G:g}
\end{align}
which admits two real solutions for~$\lambda$. Hence~$A^a$ is orthogonal to~$l_a$ which leaves two options remaining for~$A^a$ for each of the two values of~$\lambda$.

Now we want to find the rest of the eigenvectors. For that we first choose $\delta Z_1=~1$, $\delta Z_2 =~0$. Contracting the first equation with~$l_b$, and using the anti-symmetry of~$\delta F$ we get a condition for $l_a$,  \begin{align}
g_1^{ab}l_al_b = 0, \label{G:g_1}
\end{align}
which again admits two real values of~$\lambda$. Repeating the argument above, the first equation becomes,
\begin{align}
(l_al^a )A^b - (l_aA^a) l^b + g_1^{ab}l_b = 0 \label{cha_1}    
\end{align}
Since the null cones of~$g^{ab}$ and~$g_1^{ab}$ are by assumption not touching, we have $g^{ab}l_al_b\ne0$. It follows that~$A^a=-g_1^{ab}l_b/(l_cl^c)$ satisfies~\eqref{cha_1} provided that~\eqref{G:g_1} holds. Observe furthermore that~$A^a+\alpha l^a$ satisfies the same equations and results in the same Faraday tensor~$\delta F_{ab}$ for any~$\alpha$. Thus Eqn.~\eqref{G:g_1} gives two additional eigenvectors.

If we drop the assumption that the null cones of~$g^{ab}$ and~$g_1^{ab}$ are non-touching and assume that they touch at~$l_a$ then to have a solution we need that~$g^{ab}l_b$ must be proportional to~$g_1^{ab}l_b$. 

The final case is similar to the second. We choose~$\delta Z_1=0$, $\delta Z_2=1$ and  obtain
\begin{align}
 g_2^{ab}l_a l_b = 0  \label{G:g_2}    
\end{align}
and the same equations for the dual of~$\delta F_{ab}$, so we need not discuss it separately.

In summary, we have obtained the eight eigenvectors we required to satisfy the Kreiss condition and conclude that the system is strongly hyperbolic.

\subsection{Toy MHD}

Here we look at the evolution of a magnetic field~$b^a$ driven by a given velocity field~$u^a$ in a space-time~$(M,g_{ab})$. The system is
\begin{align}
    \nabla_{a}(b^{[a}u^{b]}) = 0 \label{MHD_1}
\end{align}
Here we take~$u^a$ to be time-like and of norm one, $u^au^b g_{ab} = -1$. We also take~$u^ab^b g_{ab} = 0$. This last is a gauge condition to make the solutions unique for the whole system since otherwise, if~$(u^a, b^a)$ is a solution, then~$(u^a, b^a + \eta u^a)$ also is a solution, with~$\eta$ an arbitrary function. 

We observe that there are 4 equations for 3 variables. Three of them are evolution equations for the three components of~$b^c$.  We shall see below that the other is a constraint. Thus \textbf{condition 2 } is also satisfied.

The principal part of system \eqref{MHD_1} is
\begin{align}
    \mathfrak{N}^{ba}{}_{c}\nabla_a b^c = u^{[a}\nabla_a b^{b]} = \delta^{[a}{}_c u^{b]}\nabla_a b^c.
\end{align}
It is easy to check that \textbf{condition 1} is satisfied.  The Geroch matrices are also easy to obtain as~$C^d{}_b l_d := \delta^d{}_b l_d$. They form a basis of the left kernel of~$\mathfrak{N}^{ba}{}_{c}l_a$ and, as we explained before, this means that when Eqn.~\eqref{MHD_1} is contracted with~$C^d{}_b u_d=u_b$ a constraint is generated; this is  
\begin{align}
    \nabla_a b^a - b^a a_a=0,
\end{align}
 where~$a^a \equiv u^b\nabla_b u^a$. We notice that this is the spatial divergence of~$b^a$ in disguise.

On the other hand,  the following integrability condition  $C^d{}_b \nabla_d \nabla_{a}(b^{[a}u^{b]})= \nabla_b \nabla_{a}(b^{[a}u^{b]})=0$ holds trivially, thus, the system satisfies \textbf{condition 3}.

The extended system consists of adding a term~$g_1^{ba} \nabla_a Z$ to equation~\eqref{MHD_1}, with~$g_1^{ba}$ as in the previous example and with the extra variable $Z$. Its principal part is~$u^{[a}\nabla_a b^{b]} + g_1^{bc}C^d{}_c{}{} \nabla_d Z = 0$, with~$C^a{}_b=\delta_b{}^a$. The characteristic equation is 
\begin{align}
    \frac{1}{2} (u^{a}l_a \delta b^{b}- u^{b}l_a \delta b^{a} ) + g_1^{bd} l_d \delta Z = 0 \label{MHD_cha_1}
\end{align}
where we need to solve this equation for~$l_a= -\lambda u_a+k_a$ with~$k_a$ given, and for the eigenvectors~$\delta Z$ and~$\delta b^{a}$ (with $u_a\delta b^{a}=0$).

Without loss of generality we choose $k^{a}$ such that $u^{a}k_{a}=0$, and we
rewrite the characteristic equations projecting on to~$u_{a}$ and perpendicular to it (with the projector~$h_{ab}\equiv g_{ab}+u_{a}u_{b}$). We obtain 
\begin{align*}
\frac{1}{2}k_{a}\delta b^{a}+u_{a}g_{1}^{ab}l_{b}\delta Z  & =0\\
\frac{1}{2}\lambda\delta b^{a}+h_{c}^{a}g_{1}^{cb}l_{b}\delta Z  & =0
\end{align*}
The physical solution comes from choosing~$\lambda=0$, and the eigenvectors~$\delta Z=0$ and~$\delta b^{a}$ orthogonal to $k_{a}$. Since~$\delta b^{a}$
has two possible directions we obtain two eigenvectors.

The remaining eigenvectors come from choosing~$\lambda$ such that
\begin{align}
l_{a}g_{1}^{ab}l_{b}=0, \label{MHD_g1}
\end{align}
and~$\delta Z=\frac{1}{2}\lambda$, $\delta b^{a}=h_{c}^{a}g_{1}^{cb}l_{b}$. This expression satisfies the second characteristic equation trivially, and it is easy to check
that the first one reduces to
\begin{align*}
\frac{1}{2}k_{a}\delta b^{a}+u_{a}g_{1}^{ab}l_{b}\delta Z=\frac{1}{2}%
l_{a}g_{1}^{ab}l_{b}=0. 
\end{align*}
Since, as before, there are two solution for~$\lambda$ from Eqn.~\eqref{MHD_g1}, we obtain two more eigenvectors. In summary, we have obtained the four eigenvectors we
required to satisfy the Kreiss condition and conclude that the extended system is strongly hyperbolic. Finally, we observe that we can obtain~\eqref{MHD_g1} from the integrabilty condition, multiplying Eqn.~\eqref{MHD_cha_1} by~$C^d{}_b l_d=l_b$.

\section{Conclusions}
\label{Section:Conclusions}

Similar extensions to those proposed here were previously known, starting with the divergence-cleaning used in magnetohydrodynamics and later generalized as~$\lambda$-systems for generic symmetric hyperbolic systems. To implement them it was necessary to break the covariance of the system, in the usual sense of performing a~$3+1$ decomposition. For symmetric hyperbolic systems such extensions can be obtained in our framework by committing to a frame and a reduction and adding an extra term which annihilate the time component of the constraint basis. This results in an extended symmetric hyperbolic system.

In this article we have presented an extension scheme for first order PDEs. With appropriate adaptation, however, these results can be applied to systems of order two or even higher order. We will show in future articles how to apply these ideas to gravity theories, to extend the system and to fix the gauge allowing us to reinterpret and generalize known results such as those of~\cite{Bona:2004ky, Hilditch:2013ila, Kovacs:2020ywu}.

Although the existence of a strongly hyperbolic extension is performed in Fourier space and results in system of pseudodifferential equations, our examples show that in cases of physical interest one may obtain differential extensions. These extensions furthermore retain covariance of the theory in the sense that, contrary to earlier~$\lambda$-System extensions, at least in the principal part, they do not rely on a preferred time direction, but instead the addition of other Lorentzian metric tensors. Further details and a complete proof will be given in a longer version of this work.

In our analysis we resorted to previous work to argue that the constraints, if initially satisfied, are satisfied at later times. This helped us conclude that~$Z_\Gamma$ remains zero throughout the evolution. There are however more elegant ways to show this when the constraints do not have any kernel from the left, that is, no set of zero rows in their Kronecker decomposition (see~\eqref{KD_K}). In such cases, it can be shown that the~$Z_{\Gamma}$ fields satisfy a second order evolution system which is decoupled from~$\phi^\alpha$ and has a well-posed initial value problem. Choosing these fields to vanish at the initial surface and the~$\phi^\alpha$ fields satisfy the original constraints of the system, all derivatives of~$Z_{\Gamma}$ vanish on the initial surface, in particular any transversal derivative, so the unique solution to the second order system is zero and the constraints are satisfied for all times. Unfortunately the presence of zeros may prevent the second order system from being well-posed, so that more care is needed. This will be further considered in the aforementioned longer paper.

\section*{Acknowledgements}

We thank Carlos Palenzuela for several helpful discussions associated with the article. We also thank the ``Programa de projectes de recerca amb investigadors convidats de prestigi reconegut'' of the Universitat de les Illes Balears through which Professor Oscar Reula visited the UIB, making it possible to finalise this article. OR acknowledges financial support from CONICET, SeCyT-UNC, and MinCyT-Argentina. This work was supported by the Grant PID2022-138963NB-I00 funded by MCIN/AEI/10.13039/501100011033/FEDER, UE, and by FCT Project No. UIDB/00099/2020.

\bibliographystyle{unsrt}
\bibliography{bibl_project}

\end{document}